\documentclass [10pt,a4paper,draft] {article}
\usepackage [applemac] {inputenc}
\usepackage[french] {babel}
\usepackage {amsfonts}
\usepackage {amsmath}
\usepackage {amssymb}

\begin{document}

\begin{center}
\begin{LARGE}
La propri\'et\'e de Dixmier pour les alg\`ebres de Lie de champs de vecteurs

\smallskip

\end{LARGE}
\bigskip

par : Mustapha RA\"IS  (Poitiers)
\end{center}

\bigskip
\vspace{15 mm}

\begin{large}
\noindent
\textbf{Abstract} \end{large} : Given a linear representation $\rho : \mathfrak{g} \longrightarrow
\mathfrak{g}\ell(V)$ of a Lie algebra $\mathfrak{g}$, one can define a linear representation $\rho_m :
\mathfrak{g}_m \longrightarrow \mathfrak{g}\ell(V^m)$ of the generalized Takiff algebra $\mathfrak{g}_m$. It
is proved here that the vector fields defined by $\rho_m$ on $V^m$ do have the Dixmier property if those
defined by $\rho$ have the same property. Examples where the result applies are given and in particular, those of
the adjoint or coadjoint representations of Takiff algebras.

\vskip 15mm
\begin{large}
\noindent
\textbf{1~Introduction au r\'esultat principal}
\end{large}

\medskip
                Le corps de base est $\Bbb{K} = \Bbb{R}$ ou $\Bbb{C}$. Sur les $\Bbb{K}$-espaces vectoriels, on consid\`erera
divers espaces de fonctions num\'eriques ou vectorielles, suivant le cas : l'espace des fonctions polyn\^omes
($\Bbb{K} = \Bbb{R}$ ou $\Bbb{C}$), celui des fonctions analytiques ou des fonctions $C^\infty$
($\Bbb{K} = \Bbb{R}$), celui des fonctions holomorphes enti\`eres ($\Bbb{K} = \Bbb{C}$). Ces fonctions seront dans
la suite appel\'ees les fonctions lisses, et lorsqu'il s'agira de fonctions num\'eriques, chacun des espaces
pr\'ec\'edents est une $\Bbb{K}$-alg\`ebre, qu'on notera $\mathcal{A}$.

\medskip
$\bullet$ Soient $W$ et $E$ des $\Bbb{K}$-espaces vectoriels de dimension finie. Un champ de vecteurs sur $E$,
avec param\`etres dans $W$, est une fonction lisse $a : W \times E \longrightarrow E$, ou un op\'erateur
diff\'erentiel lin\'eaire du premier ordre sur $E$ :
$$
(L\varphi)(w,e) = \bigg({d\over dt}\bigg)_0\ \varphi(e+ta(w,e))
$$
o\`u $\varphi$ est une fonction num\'erique lisse sur $E$.

\medskip
$\bullet$ Consid\'erons en particulier une repr\'esentation lin\'eaire $\rho : \mathfrak{h} \longrightarrow
\mathfrak{g}\ell(E)$ d'une alg\`ebre de Lie $\mathfrak{h}$, de dimension finie. Chaque \'el\'ement $x$ de
$\mathfrak{h}$ d\'efinit un champ de vecteurs $a_x : E \longrightarrow E,\ a_x(e) = \rho(x).e$, i.e. un
op\'erateur diff\'erentiel $L_x$, avec $(L_x\varphi)(e) = \big({d\over dt}\big)_0\ \varphi(e+t\rho(x).e)$ pour
toute fonction lisse sur $E$. L'ensemble $\{L_x| x\in \mathfrak{h}\}$ de ces champs de vecteurs (les champs de
Killing de $\rho$ selon la terminologie de \cite{D} (\no 19.3)) est une alg\`ebre de Lie de champs de vecteurs sur
$E$ et l'ensemble des fonctions lisses $\varphi : E \longrightarrow \Bbb{K}$ telles que : $L_x\varphi = 0$ pour tout
$x$, est la sous-alg\`ebre $\mathcal{A}^\mathfrak{g}$ des fonctions lisses $\mathfrak{g}$-invariantes.

\vskip 7mm
\noindent
\textbf{D\'efinition} :  On dira que la repr\'esentation $(\mathfrak{h}, \rho, E)$ a la
\underline{propri\'et\'e de Dixmier} lorsque :

        Pour tout espace de param\`etres $W$, pour tout champ de vecteurs avec param\`etres $a : W \times E
\longrightarrow E$ qui annule toutes les fonctions $\mathfrak{g}$-invariantes sur $E$,

\medskip
$\bullet$ il existe une fonction lisse
$b : W \times E \longrightarrow \mathfrak{h}$ telle que :
$$
        a(w,e) = \rho(b(w,e)).e \quad ((w,e) \in W \times E)
$$

\medskip
$\bullet$ ou encore, de fa\c con \'equivalente, il existe des fonctions lisses $\psi_j$ sur $W \times E$ telles que :
$$
                a(w,e) = \sum_j \psi_j \ (w,e) \rho(x_j).e
$$

\medskip
o\`u $(x_j)_j$ est une base de $\mathfrak{h}$.

\medskip
$\bullet$ ou encore, de fa\c con \'equivalente : $a$ est combinaison lin\'eaire, \`a coefficients dans l'alg\`ebre des
fonctions lisses sur $W \times E$, des champs de Killing de la repr\'esentation $\rho$.

\bigskip
- Soient  $\mathfrak{g}$ une $\Bbb{K}$-alg\`ebre de Lie de dimension finie, $m$ un entier $\geq 0$, et
$\mathfrak{g}_m$ l'alg\`ebre de Lie construite dans \cite{R-T} (et connue depuis sous l'appellation d'alg\`ebre
de Takiff g\'en\'eralis\'ee). Les \'el\'ements de $\mathfrak{g}_m$, d\'esign\'es par des lettres capitales,
s'\'ecrivent :
$$
        X = x_0 + x_1 T+\cdots + x_m T^m = (x_0, x_1,\ldots ,x_m)
$$
o\`u $T$ est une ind\'etermin\'ee, et les $x_j$ sont des \'el\'ements de $\mathfrak{g} = \mathfrak{g}_0$. Le
crochet de Lie de $\mathfrak{g}_m$ est d\'efini par :
$$
        \big[\sum \ x_r\, T^r, \sum \ y_s\, T^s\big] = \sum_{r,s}\, \big[x_r, y_s\big]T^{r+s}\ mod\  T^{m+1}
$$

\bigskip
- Soit maintenant $\rho : \mathfrak{g} \longrightarrow \mathfrak{g}\ell(V)$ une repr\'esentation lin\'eaire de
$\mathfrak{g}$ dans un espace vectoriel $V$. On d\'efinit l'espace vectoriel $V_m$ comme \'etant l'espace des :
$$
        F = f_0 + f_1\ T+\cdots + f_mT^m = (f_0,f_1,\ldots ,f_m)
$$
o\`u les $f_j$ sont des \'el\'ements de $V$, et une repr\'esentation lin\'eaire $\rho_m$ de $\mathfrak{g}_m$
dans $V_m$ :
$$
        \rho_m(x_rT^r).f_sT^s = (\rho(x_r).f_s)T^{r+s} \ mod\  T^{m+1}
$$

\noindent
Le r\'esultat principal obtenu ici est le

\vskip 7mm
\noindent
\textbf{Th\'eor\`eme} :  \textit {Lorsque la repr\'esentation $(\mathfrak{g},\rho,V)$ a la propri\'et\'e de
Dixmier, il en est de m\^eme de la repr\'esentation $(\mathfrak{g}_m,\rho_m,V_m)$.}


\vskip 14mm
\begin{large}
\noindent
\textbf{2~Construction de fonctions lisses $\mathfrak{g}_m$-invariantes sur $V_m$}
\end{large}

\vskip 7mm
\noindent
\textbf{2.1. } Une fonction lisse $\Phi : V_m \longrightarrow \Bbb{K}$ est
$\mathfrak{g}_m$-invariante si et seulement si :
$$
        \bigg({d\over dt}\bigg)_0\Phi(F+t\rho_m(X).F) = 0
$$
pour tous $\displaystyle F = \sum_r\, f_r\, T^r$ et $X = \sum_s\, x_s\, T^s$, c'est-ˆ-dire si et seulement si :
$$
        \bigg({d\over dt}\bigg)_0 \Phi(f_0 + t\rho(x_0)f_0, f_1 + t\{\rho(x_0)f_1+\rho(x_1)f_0\},\ldots ,
$$
$$
        f_m+t\{\rho(x_0)f_m+\rho(x_1)f_{m-1}+\cdots + \rho(x_m)f_0\}) = 0
$$

Supposons en particulier que la fonction $\Phi$ ne d\'epende pas de la composante $f_m$, i.e. qu'il existe une
fonction lisse $\theta$ sur $V_{m-1}$ telle que : $\Phi(f_0,f_1,\ldots , f_m) = \theta(f_0, f_1,\ldots , f_{m-1})$.
La formule ci-dessus exprime l'invariance de $\Phi$ sous l'action de $\mathfrak{g}_m$ de la mani\`ere suivante
:
$$
        \bigg({d\over dt}\bigg)_0 \ \theta(f_0 + t\rho(x_0).f_0,\ldots ,f_{m-1} + t\{\rho(x_0)f_{m-1} + \rho(x_1)f_{m-2}
+\cdots +
\rho(x_{m-1}).f_0\}) = 0
$$

On notera donc, pour un usage ult\'erieur, la

\vskip 5mm
\noindent
\underline{Remarque} : La fonction $\Phi$ est $\mathfrak{g}_m$-invariante sur $V_m$ si et seulement si la
fonction $\theta$ est $\mathfrak{g}_{m-1}$-invariante sur $V_{m-1}$.

\vskip 7mm
\noindent
\textbf{2.2.}~Pour construire des fonctions lisses $\mathfrak{g}_m$-invariantes sur $V_m$, on adapte le
proc\'ed\'e utilis\'e pour les fonctions \underline{polyn\^omes} dans \cite{R-T} (Lemmes 3.2 et 3.5), aux
fonctions lisses non n\'ecessairement polynomiales, en effectuant quelques modifications qui semblent
n\'ecessaires. Pour simplifier les notations, et pour autant que cela n'introduise pas d'ambigu\"it\'e, on \'ecrira :
$$
\begin{array}{llll}
        \rho_m(X)F &= \rho(X)F &= X.F         &(X\in \mathfrak{g}_m,\ F \in V_m)\\
\smallskip\\
\rho(x)f        &=x.f        &        &(x \in \mathfrak{g},\ f \in V)\\
\end{array}
$$

\vskip 7mm
\noindent
\textbf{2.3.}~Soient $F = (f_0, f_1,\ldots , f_m)$ dans $V_m,\ g : \Bbb{K} \longrightarrow V$ la fonction
d\'efinie par : $g(t) = f_0 + tf_1+\cdots+ t\  {}^m\!f_m$, et $\varphi : V \longrightarrow \Bbb{K}$ une fonction
lisse. Pour chaque entier naturel $k$, on pose :
$$
        \Phi_k(f_0,\ldots , f_m) = {1\over k!}\ \bigg({d\over dt}\bigg)^k_0\  \varphi(g(t))
$$
et on obtient ainsi des fonctions lisses $\Phi_k : V_m \longrightarrow \Bbb{K}$

\vskip 7mm
\noindent
\textbf{2.4.} \textbf{Proposition} :  \textit{Soit $k$ un entier tel que $0 \leq k \leq m$}

\smallskip
\begin{enumerate}
\item \textit{ $\Phi_k$ ne d\'epend que des composantes $f_0,f_1,\ldots , f_k$, et pr\'ecis\'ement :
$$
        \Phi_k(F) = {1\over k!} \bigg({d\over dt}\bigg)^k_0\ \varphi(f_0+tf_1+\cdots + t \ {}^k\!f_k)
$$}

\medskip
\item \textit{$\Phi_k$ est une fonction lin\'eaire de la composante $f_k$, et pr\'ecis\'ement :
$$
        \Phi_k(F) = \ <d\varphi(f_0),f_k> + \ \psi_k(f_0,f_1,\ldots , f_{k-1})\qquad (*)
$$
o\`u $\psi_k : V_{k-1}\longrightarrow \Bbb{K}$ est une fonction lisse sur $V_{k-1}$.}

\medskip
\item \textit{Si $\varphi$ est une fonction lisse $\mathfrak{g}$-invariante sur $V$, la fonction $\Phi_k : V_m
\longrightarrow \Bbb{K}$ est lisse et $\mathfrak{g}_m$-invariante sur $V_m$.}
\end{enumerate}

\bigskip
\noindent
\textsc{D\'emonstration} :

$\bullet$ D'apr\`es la formule classique de Faa-di-Bruno
$$
        {1\over k!} \big({d\over dt}\big)^k\ \varphi o g = \sum {1\over q_1!q_2!\cdots q_k!} (d^q\varphi)og.({g'\over
1!})^{[q_1]}\  ({g''\over 2!})^{[q_2]}\cdots ({g^{(k)}\over k!})^{[q_k]})
$$
o\`u $q = q_1+q_2+\cdots + q_k$ et la somme est prise sur tous les entiers $q_1,q_2,\ldots ,q_k$ v\'erifiant :
$q_1 + 2q_2 +\cdots + kq_k = k$ (les notations \'etant celles habituelles du calcul diff\'erentiel). Par exemple
$(d^q\varphi)(g(t))$ est une
$q$-forme multilin\'eaire sym\'etrique, qu'on \'evalue sur $q$ vecteurs o\`u ${g'(t)\over 1!}$ figure $q_1$ fois,
${g''(t)\over 2!}$ figure $q_2$ fois,$\ldots , {g^{(k)}(t)\over k!}$ figure $q_k$ fois.

\smallskip
On a : $\displaystyle {1\over r!}\ g^{(r)}(0) = f_r$ et :
$$
        \Phi_k (f_0,f_1,\ldots , f_m) = \sum {1\over q_1!\cdots q_k!}\ (d^q\varphi)(f_0).(f_0^{[q_1]},\ldots
,f_k^{[q_k]})
$$
ne d\'epend que des composantes $f_0,f_1,\ldots , f_k$. Ceci d\'emontre la partie 1.

\vskip 5mm
$\bullet$ Dans la somme ci-dessus, il y a le terme correspondant \`a $q_1 = \cdots = q_{k-1} = 0$ et $q_k = 1$,
et dans tous les autres termes, on a : $q_1 + 2q_2 +\cdots + (k-1)q_{k-1} = k$ et $q_k = 0$. Ainsi :
$$
        \Phi_k(f_0,f_1,\ldots , f_m) =\, <d\varphi(f_0),f_k>\ +\  \psi_k(f_0,f_1,\ldots , f_{k-1})
$$
o\`u $\psi_k$ est une fonction lisse sur $V_{k-1}$. Ceci d\'emontre la partie 2.

\vskip 5mm
$\bullet$ Soit $\varphi$ une fonction lisse $\mathfrak{g}$-invariante sur $V$. Pour tous $t$ et $s$ dans
$\Bbb{K}$, $x_0,x_1,\ldots , x_m$ dans $\mathfrak{g}$, et $f_0,f_1,\ldots ,f_m$ dans $V$, on a :
$$
        \bigg({d\over ds}\bigg)_0\ \varphi(\sum t^j f_j + s\rho(\sum\, t^i x_i).\sum\, t^j f_j) = 0.
$$
On a : $\displaystyle \rho(\sum_{0\leq i \leq m}\ t^i x_i).\sum_{0\leq j \leq m}\ t^j f_j = \sum_{0\leq j \leq m}\
t^j h_j + t^{m+1}k(t)$, avec :
$$
        h_j = \sum_{0 \leq r \leq j}\ x_r.f_{j-r}.
$$
D'o\`u :
$$
        \sum t^jf_j + s\, \rho(\sum t^ix_i).\sum t^jf_j = \sum t^j(f_j + sh_j) + st^{m+1}k(t)
$$
et lorsque $0 \leq k \leq m$ :
$$
        \bigg({d\over dt}\bigg)^k_0\ \varphi (\sum t^jf_j + s \rho(\sum t^ix_i).\sum t^jf_j)
$$

$$
= \bigg({d\over dt}\bigg)^k_0\ \varphi(\sum t^j(f_j + sh_j) = \Phi_k(f_0 + sh_0 , f_1 + sh_1,\ldots , f_m + sh_m))
$$
On a donc :
$$
        \big({d\over ds}\big)_0 \, \Phi_k(f_0 + sh_0, f_1 + sh_1,\ldots , f_m + sh_m) = \big({d\over dt}\big)^k_0\
({d\over ds})_0\
\varphi(\sum t^jf_j + s\rho (\sum t^ix_i).\sum t^jf_j) = 0
$$
Comme : $(f_0 + sh_0, f_1 + sh_1,\ldots , f_m + sh_m) = (f_0,\ldots , f_m) + s \rho_m(x_0,\ldots ,
x_m).(f_0,\ldots , f_m)$. On voit que les $\Phi_k$ sont $\mathfrak{g}_m$-invariantes.

\vskip 14mm
\begin{large}
\noindent
\textbf{3~La propri\'et\'e de Dixmier pour la repr\'esentation $(\mathfrak{g}_m, \rho_m, V_m)$}
\end{large}

\medskip
        Il s'agit de montrer que $(\mathfrak{g}_m, \rho_m, V_m)$ a la propri\'et\'e de Dixmier, sous l'hypoth\`ese que
$(\mathfrak{g}, \rho, V)$ a cette m\^eme propri\'et\'e. La d\'emonstration se fait par r\'ecurrence sur $m$,
sachant que le cas $m = 0$ correspond \`a l'hypoth\`ese. On suppose $m > 0$.

\medskip
$\bullet$ Soit $a = (a_0, a_1,\ldots , a_m) : W \times V_m \longrightarrow  V_m$ un champ de vecteurs qui
annule toutes les fonctions $\mathfrak{g}_m$-invariantes sur $V_m$, {\it i.e.} tel que :

$$
        \bigg({d\over ds}\bigg)_0 \Phi(f_0 + sa_0(w,f_0,f_1,\ldots ,f_m)),\ldots , f_m + sa_m(w,f_0,f_1,\ldots ,f_m)) = 0
$$
pour toute fonction $\psi : V_m \longrightarrow \mathbb{K}$, $\mathfrak{g}_m$-invariante.

        On applique ceci, en particulier, \`a toutes les fonctions $\mathfrak{g}_m$-invariantes :
$$
        \Phi(f_0,\ldots ,f_m) = \theta(f_0,\ldots , f_{m-1})
$$
qui ne d\'ependent que des composantes $f_0,\ldots , f_{m-1}$. D'apr\`es la remarque 2.1, il vient :
$$
        \bigg({d\over ds}\bigg)_0 \theta(f_0 + sa_0(w,\ldots , f_m),\ldots , f_{m-1} + sa_{m-1}(w, f_0,\ldots ,f_m)) = 0
$$
pour toute fonction $\theta$, $\mathfrak{g}_{m-1}$-invariante sur $V_{m-1}$. Par l'hypoth\`ese de
r\'ecurrence, il existe des fonctions lisses $b_j : W \times V_m \longrightarrow \mathfrak{g}\ (0 \leq j \leq
m-1)$ telles que :
$$
        (a_0, a_1,\ldots , a_{m-1}) = \rho_{m-1}(b_0, b_1,\ldots , b_{m-1}).(f_0,\ldots , f_{m-1}).
$$

        On v\'erifie imm\'ediatement que :
$$
        \rho_m(b_0 + b_1 T +\cdots + b_{m-1} T^{m-1}).(f_0 + f_1 T +\cdots + f_m T^m) = (a_0, a_1, \ldots ,
a_{m-1},0)+(0,0,\ldots , 0,c_m)
$$
avec $c_m = b_0.f_m + b_1.f_{m-1} +\cdots + b_{m-1}.f_1$. D'o\`u :
$$
        (a_0,a_1,\ldots , a_m) = (a_0,a_1,\ldots , a_{m-1},0)+(0,0,\ldots , 0, a_m)
$$
$$
= \rho_m(b_0,\ldots , b_{m-1},
0).(f_0,f_1,\ldots , f_m)+ (0,0,\ldots , 0, a_m - c_m)
$$
et le champ de vecteurs $(0,0,\ldots , 0, a_m-c_m)$ annule les fonctions $\mathfrak{g}_m$-invariantes. En
particulier, pour toute fonction lisse $\mathfrak{g}$-invariante $\varphi : V \longrightarrow \mathbb{K}$, on a
:
$$
        ({d\over ds})_0 \, \Phi_m(f_0,f_1,\ldots , f_{m-1}, f_m + s(a_m - c_m)) = 0
$$
avec $\Phi_m(f_0,\ldots , f_{m-1}, f_m) = \, <d\varphi(f_0),f_m> + \, \psi_m(f_0,f_1,\ldots , f_{m-1})$  (on
applique la formule ($\star$) de 2.4). Donc :
$$
        <d\varphi(f_0), a_m-c_m>\, = 0
$$
et \`a nouveau la propri\'et\'e de Dixmier de $(\mathfrak{g}, \rho, V)$ assure l'existence d'une fonction lisse
$b_m : W  \times V_m \longrightarrow \mathfrak{g}$ telle que
$$
        a_m = c_m + b_m.f_0 = b_0.f_m + b_1.f_{m-1} + \cdots + b_{m-1}.f_1 + b_m.f_0
$$
Ainsi : $(a_0,a_1,\ldots , a_m) = \rho_m(b_0,\ldots , b_m).(f_0,f_1,\ldots , f_m)$, cqfd.

\vfill\eject
\begin{large}
\noindent
\textbf{4~\underline{Notes bibliographiques}}
\end{large}

\bigskip
\noindent
\textbf{4.1.}~A l'origine (cf. \cite{Dix}), J. Dixmier examine la repr\'esentation adjointe $(\mathfrak{g}, ad,
\mathfrak{g})$ d'une alg\`ebre de Lie semi-simple sur $\mathbb{K} = \mathbb{R}$ ou $\mathbb{C}$ et
d\'emontre la ``propri\'et\'e de Dixmier'' (toutefois pour les champs de vecteurs sans param\`etre), pour
toutes les alg\`ebres $\mathcal{A}$ de l'introduction \`a l'exception, lorsque $\mathbb{K} = \mathbb{R}$, de
celle des fonctions de classe $C^\infty$. En fait l'extension aux champs de vecteurs avec param\`etres est
imm\'ediate, et en particulier, on peut supposer $\mathfrak{g}$ r\'eductive ci-dessus.

\vskip 5mm
\noindent
\textbf{4.2.}~La d\'emonstration de la propri\'et\'e de Dixmier pour la repr\'esentation adjointe
$(\mathfrak{g}_m, ad_m, \mathfrak{g}_m)$ de l'alg\`ebre de Takiff g\'en\'eralis\'ee, dans le cas o\`u
$\mathfrak{g}$ est semi-simple et $\mathcal{A}$ est l'alg\`ebre des fonctions polyn\^omes, se trouve dans
\cite{Ra}.

\vskip 5mm
\noindent
\textbf{4.3.}~Suite \`a une question que je lui ai pos\'ee, A. Bouaziz a d\'emontr\'e la propri\'et\'e de Dixmier
pour la repr\'esentation adjointe d'une alg\`ebre de Lie r\'eductive r\'eelle et pour l'alg\`ebre $\mathcal{A}$
des fonctions de classe $C^\infty$ (cf. \cite{Bou}).

        Compte tenu du th\'eor\`eme d\'emontr\'e ici, la repr\'esentation adjointe d'une alg\`ebre de Takiff
$\mathfrak{g}_m$, construite sur une alg\`ebre de Lie r\'eductive $\mathfrak{g}$, a la propri\'et\'e de
Dixmier pour toutes les alg\`ebres $\mathcal{A}$ de fonctions lisses.

\vskip 5mm
\noindent
\textbf{4.4.}~La propri\'et\'e de Dixmier a \'et\'e d\'emontr\'ee par D. Panyushev (\cite{Pa}) dans le cas o\`u
$\mathcal{A}$ est l'alg\`ebre des fonctions polyn\^omes, pour les alg\`ebres de Takiff $\mathfrak{g}_m$,
sous la condition que $\mathfrak{g}$ appartienne \`a une classe particuli\`ere d'alg\`ebres de Lie, celle dite
des ``3-wonderful Lie algebras''. Panyushev d\'emontre que $\mathfrak{g}_m$ est une  ``3-wonderful Lie
algebra'' d\`es que $\mathfrak{g}$ elle-m\^eme appartient \`a cette m\^eme classe d'alg\`ebres de Lie, et
qu'en particulier les alg\`ebres de Lie r\'eductives complexes sont des  ``3-wonderful Lie algebras''.

\vskip 14mm
\begin{large}
\noindent
\textbf{5~\underline{Quelques remarques}}
\end{large}

\bigskip
\noindent
\textbf{5.1.}~Sur  $V_m$, on a construit, pour chaque $\varphi$ dans $\mathcal{A}^{\mathfrak{g}}$, $(m+1)$
fonctions $\mathfrak{g}_m$-invariantes $\Phi_0, $ $\Phi_1,\ldots , \Phi_m$. Lorsque $\varphi$ d\'ecrit
$\mathcal{A}^{\mathfrak{g}}$, l'ensemble de ces fonctions est en fait une sous-alg\`ebre $\mathcal{B}$ de
$\mathcal{A}^{\mathfrak{g}_m}(V_m)$, et un examen attentif de la d\'emonstration donn\'ee dans \textbf{3}
ci-dessus montre que pour qu'un champ de vecteurs annule toutes les fonctions $\mathfrak{g}_m$-invariantes, il
suffit qu'il annule les fonctions appartenant \`a $\mathcal{B}$. En un certain sens \`a pr\'eciser, la sous-alg\`ebre
$\mathcal{B}$ est dense dans $\mathcal{A}$.

\vskip 7mm
\noindent
\textbf{5.2.}~Soit $\rho : \mathfrak{g} \longrightarrow \mathfrak{g}\ell(V)$ une repr\'esentation de
$\mathfrak{g}$ dans $V$. Un champ de vecteurs $a : W \times V \longrightarrow  V$ est dit tangent aux orbites
de $\mathfrak{g}$, lorsque :
$$
        a(w,v) \in \rho(\mathfrak{g}).v \quad \textrm {pour presque tout} \ v \ \textrm {dans} \ V
$$
(presque tout \'etant relatif suivant le cas \`a la topologie de Zariski ou \`a la topologie num\'erique).

        Supposons que $(\mathfrak{g}, \rho, V)$ ait la propri\'et\'e de Dixmier. Il est clair alors qu'un champ de
vecteurs sur $V$ annule les fonctions $\mathfrak{g}$-invariantes si et seulement s'il est tangent aux
$\mathfrak{g}$-orbites.

\vskip 7mm
\noindent
\textbf{5.3.}~Soit  $\rho : \mathfrak{g} \longrightarrow \mathfrak{g}\ell(V)$ une repr\'esentation lin\'eaire de
$\mathfrak{g}$ et soit $\theta : V \longrightarrow V'$ une bijection lin\'eaire de $V$ sur un espace vectoriel
$V'$. La formule :
$$
        \tau(x) = \theta \circ \rho(x) \circ \theta^{-1}\quad (x \in \mathfrak{g})
$$
d\'efinit une repr\'esentation $\tau$ de $\mathfrak{g}$ dans $V'$, (\'equivalente \`a $\rho$). On v\'erifie
ais\'ement que $(\mathfrak{g}, \rho, V)$ a la propri\'et\'e de Dixmier si et seulement si $(\mathfrak{g}, \tau,
V')$ a cette m\^eme propri\'et\'e. Une application de cette remarque est faite ci-dessous.

\vskip 7mm
\noindent
\textbf{5.4.}~Soit $\mathfrak{g}$ une alg\`ebre de Lie quadratique, de sorte que l'alg\`ebre de Takiff
$\mathfrak{g}_m$ est elle-m\^eme quadratique (voir \cite {R-T}, 3.8). Autrement dit les repr\'esentations
adjointe et coadjointe de $\mathfrak{g}_m$ sont \'equivalentes. Par suite, lorsque $(\mathfrak{g}, ad,
\mathfrak{g})$ a la propri\'et\'e de Dixmier, la repr\'esentation coadjointe de $\mathfrak{g}_m$ a cette
m\^eme propri\'et\'e (et vice-versa).

\vskip 7mm
\noindent
\textbf{5.5.}~Soit $\mathfrak{g}$ une alg\`ebre de Lie, dont la repr\'esentation coadjointe a la propri\'et\'e de
Dixmier. D'apr\`es le th\'eor\`eme principal d\'emontr\'e ici, la repr\'esentation $(\mathfrak{g}_m,
(ad^*)_m, (\mathfrak{g}^*)_m)$ a cette m\^eme propri\'et\'e. Mais cette repr\'esentation de
$\mathfrak{g}_m$ dans $(\mathfrak{g}^*)_m$ (identifi\'e de mani\`ere \'evidente \`a $(\mathfrak{g}_m)^*$)
\underbar{n'est pas} la repr\'esentation coadjointe de $\mathfrak{g}_m$. Toutefois, notons $\rho$ la
repr\'esentation $(\mathfrak{g}_m, (ad^*)_m, (\mathfrak{g}^*)_m)$ et $\tau$ la repr\'esentation
$(\mathfrak{g}_m, (ad_m)^*, \mathfrak{g}^*_m)$. Avec $X = (x_0,\ldots ,x_m)$ dans $\mathfrak{g}_m$ et $F
= (f_0, f_1,\ldots ,f_m)$ on a :

        $\rho(X).F = (g_0, g_1, \ldots , g_m)$ avec :
$$
        g_k = ad^*(x_0).f_k + ad^*(x_1)f_{k-1} + \cdots + ad^*(x_k).f_0
$$
tandis que : (formules extraites de (\cite {R-T}, 1.5) :

        $\tau(x).G = (g'_0, g'_1,\ldots , g'_m)$ avec
$$
        g'_k = ad^*(x_0).f_k + ad^*(x_1).f_{k+1} +\cdots + ad^*(x_{m-k}).f_m
$$
On v\'erifie alors que :
$$
        \tau(X) = \theta \circ \rho(X)\circ \theta
$$
o\`u $\theta(f_0,\ldots , f_m) = (f_m, f_{m-1},\ldots , f_0)$.

\medskip
Ainsi : lorsque la repr\'esentation coadjointe de $\mathfrak{g}$ a la propri\'et\'e de Dixmier, la
repr\'esentation coadjointe de $\mathfrak{g}_m$ a cette m\^eme propri\'et\'e.

\vskip 14mm
\begin{large}
\noindent
\textbf{6~\underline{Autres exemples}}
\end{large}

\bigskip
\noindent
\textbf{6.1.}~Dans le  cas de l'alg\`ebre $\mathcal{A}$ des fonctions polyn\^omes, l'\'etude de la propri\'et\'e
de Dixmier pour des repr\'esentations $(\mathfrak{g}, \rho, V)$, o\`u $\rho$ n'est ni la repr\'esentation
adjointe ni la repr\'esentation coadjointe d'une alg\`ebre de Lie, a fait l'objet de travaux de Levasseur et
Stafford, Levasseur et Ushirobira, et plus r\'ecemment de Panyushev. On pourra consulter \cite {Pa} et sa
bibliographie.

\vskip 7mm
\noindent
\textbf{6.2.}~Dans le cas des alg\`ebres de fonctions $C^\infty$, des r\'esultats n\'egatifs sont d\'ecrits dans
(\cite {Ra}, 3) et au moins un r\'esultat positif qui est celui o\`u $\mathfrak{g} = \mathbb{R}$ est ``le'' radical
nilpotent de $sl(2)$, et $\rho$ est la restriction \`a $\mathfrak{g}$ de la repr\'esentation irr\'eductible de
dimension $\geq 3$ de $sl(2)$ (\cite {Ra}, 3.2.2). En appliquant le th\'eor\`eme principal, on peut multiplier
les exemples de repr\'esentations des alg\`ebres $\mathbb{R}^k$  ayant la propri\'et\'e de Dixmier.

\vskip 7mm
\noindent
\textbf{6.3.}~On traite ici le cas d'un espace sym?trique de rang 1. Soient $\mathfrak{g} = so(n)$ et $\rho :
\mathfrak{g} \longrightarrow \mathfrak{g}\ell(\mathbb{R}^n)$ la repr\'esentation naturelle de $so(n)$ dans
$\mathbb{R}^n$. Soit $L = \sum a_i {\partial\over \partial x_i}$ un champ de vecteurs sur $\mathbb{R}^n$,
\`a coefficients $C^\infty$, qui annule la fonction $\mathfrak{g}$-invariante :
$$
        Q(x_1,\ldots , x_n) = {1\over 2} (x^2_1 +\cdots + x^2_n).
$$
On a donc $\sum a_i x_i = 0$. On pose :
$$
        \omega = dQ = \sum x_i dx_i
$$
$$
        \alpha = \sum_i (-1)^{i-1} a_i dx_1 \wedge\cdots \wedge \big(dx_i\big)^\vee \wedge\cdots \wedge dx_n
$$
de sorte que $\displaystyle \omega_\wedge \alpha = \big(\sum_i  a_ix_i\big)dx_1 \wedge \cdots \wedge
dx_n$ et $LQ = 0$ si et seulement si $\omega_\wedge \alpha = 0$. D'apr\`es \cite {M}, il existe une
$(n-2)$-forme diff\'erentielle $\beta$ sur $\mathbb{R}^n$, \`a coefficients $C^\infty$, telle que : $\alpha =
\omega \wedge \beta$. Il en r\'esulte qu'il existe une matrice antisym\'etrique $(b_{ij})$ \`a coefficients
$b_{ij}$ fonctions $C^\infty$ sur $\mathbb{R}^n$ telle que : $\displaystyle a_i = \sum_j\, b_{ij} x_j\ (1\leq i \leq
n)$. On a donc :
$$
        \sum_i a_i {\partial\over \partial x_i} = \sum_{i,j}\, b_{ij} x_j {\partial\over \partial x_i} = \sum_{i<j}
b_{ij}(x_j {\partial\over \partial x_i} - x_i {\partial\over \partial x_j})
$$
et $(so(n), \rho, \mathbb{R}^n)$ a la propri\'et\'e de Dixmier.

\vskip 15mm

\renewcommand{\refname}{Bibliographie}     


\begin{thebibliography}{aaa}


        \bibitem[Bou]{Bou} \textsc{BOUAZIZ A.},  \emph{Sur les distributions covariantes dans les alg\`ebres de Lie
r\'eductives}. J. Funct. Analysis, vol. 257, \no 10 (2009), 3203-3217.

        \bibitem[D]{D} \textsc{DIEUDONN\'E J.},  \emph{El\'ements d'analyse}. Tome IV, Gauthier-Villars, Paris, (1971).

        \bibitem[Dix]{Dix} \textsc{DIXMIER J.},  \emph{Champs de vecteurs adjoints sur les groupes et alg\`ebres de Lie
semi-simples}. J. Reine Angew. Math. 309 (1979), 183-190.

        \bibitem[M]{M} \textsc{MOUSSU R.},  \emph{Le th\'eor\`eme de de Rham sur la division des formes}.
C.R. Acad. Sci. Paris, 280, (1975), \no 6, 329-332.

        \bibitem[Pa]{Pa} \textsc{PANYUSHEV D.},  \emph{Adjoint vector fields and differential operators on
representation spaces}. Bull. Lond. Math. Soc. 40 (2008), \no 6, 1045-1064.

        \bibitem[R-T]{R-T} \textsc{RA\"{I}S M. \& TAUVEL P.},  \emph{Indice et polyn\^omes invariants pour certaines
alg\`ebres de Lie}. J. Reine Angew. Math. 425 (1992), 123-140.

        \bibitem[Ra]{Ra} \textsc{RA\"{I}S M.},  \emph{Notes sur la notion d'invariant caract\'eristique}. Arxiv :
0707.0782.



\end{thebibliography}
\end {document}